\documentclass[12pt]{amsart}
\usepackage{amsmath}
\usepackage{amssymb}
\usepackage{amscd}
\usepackage[mathscr]{eucal}
\usepackage{latexsym}

\newtheorem{thm}{Theorem}[section]

\newtheorem{prop}[thm]{Proposition}

\theoremstyle{definition}

\newtheorem{rem}[thm]{Remark}

\newtheorem{defn}[thm]{Definition}

\newtheorem{ex}[thm]{Example}

\theoremstyle{remark}

\numberwithin{equation}{section}

\def\-rig{\text{\rm -rig}}
\def\-log{\text{\rm -log}}
\def\-dif{\text{\rm -dif}}

\def\ord{\text{\rm ord}}

\def\Gal{\text{\rm Gal}}

\def\cos{\text{\rm cos}}

\def\Cl{\text{\rm Cl}}
\def\PSL{\text{\rm PSL}}

\def\Spec{\text{\rm Spec}\,}

\def\lim{\text{\rm lim}}

\begin{document} 

\title[On the topological aspects of arithmetic elliptic curves]
{On the topological aspects of arithmetic elliptic curves}

\author[Kazuma Morita]{Kazuma Morita}
\address{Department of Mathematics, Hokkaido University, Sapporo 060-0810, Japan}
\email{morita@math.sci.hokudai.ac.jp}

\date{\today}

\maketitle
{\bf Abstract.}
In this short note, we shall construct a certain topological family which contains all elliptic curves over $\mathbb{Q}$ and, as an application, show that this family provides 
some geometric interpretations of the Hasse-Weil L-function of an elliptic curve over $\mathbb{Q}$ whose Mordell-Weil  group is of rank $\leq1$.
\section{Introduction} 
For any elliptic curve $E$ over $\mathbb{Q}$, there exists a rational newform $f$ such that we have $L(E,s)=L(f,s)$ and, in particular, the Fourier expansion of $f$ tells us the eigenvalues of the Frobenius operator acting on the Tate module of  the strong Weil curve modulo $p$. In this paper, we shall deform the Fourier expansion of $f$ with respect to the arguments $\verb+{+\theta_{p}\verb+}+_{p}$ of these eigenvalues and construct a topological family attached to these deformed differential forms. This family contains all elliptic curves over $\mathbb{Q}$ up to  isogeny  and we expect that we can deduce the arithmetic facts by  using the topological methods. Actually, as an application, if $E$ is an  elliptic curve over $\mathbb{Q}$ whose Mordell-Weil  group is of rank $\leq 1$, we will show  that this family provides 
some geometric interpretations of the Hasse-Weil L-function of $E$.

{\bf Acknowledgments.}
The author would like to thank Professor Masanori Asakura and Iku Nakamura for useful discussions. This research was partially supported by JSPS Grant-in-Aid for Research Activity Start-up. 
\section{Review of the classical theory}
Let $\mathbb{H}$ be the upper half-plane and $\mathbb{H}^{*}=\mathbb{H}\cup\mathbb{Q}\cup \verb+{+\infty\verb+}+$ be the extended upper half-plane which is obtained by adding the cusps $\mathbb{Q}\cup \verb+{+\infty\verb+}+$.  The modular group $\Gamma=\PSL_{2}(\mathbb{Z})$  acts discontinuously on $\mathbb{H}$ via linear fractional transformations. Let $\Gamma_{0}(N)$ denote the congruence subgroup $$\bigl\{
\begin{pmatrix}
a & b\\
c & d\\
\end{pmatrix}\in \Gamma\mid c\equiv  0 \ (\text{mod} \ N)\bigr\}$$
of $\Gamma$. The space of cusp forms of weight $2$ for $\Gamma_{0}(N)$ will be denoted by $S_{2}(N)$. Then, every cusp form $f(z)\in S_{2}(N)$ ($z\in \mathbb{H}$) has the Fourier expansion 
$$f(z)=\sum_{n=1}^{\infty}a_{n}(f)q^{n}\qquad (a_{n}(f)\in\mathbb{C}, \ q=e^{2\pi iz}).$$ 
We say that $f(z)$ is a normalized cusp form if we have $a_{1}(f)=1$. On the other hand, the space of cusp forms $S_{2}(N)$ is equipped with the  Hecke operators: 
\begin{itemize}
\item 
$T_{p}: f(z) \mapsto pf(pz)+\frac{1}{p}\sum_{r=0}^{p-1}f(\frac{z+r}{p})$ $\quad (p\nmid N$ ($p$: prime))
\item 
$U_{p}: f(z) \mapsto \frac{1}{p}\sum_{r=0}^{p-1}f(\frac{z+r}{p})$ $\mspace{48mu}$ $\qquad (p\mid N$ ($p$: prime)).
\end{itemize}
Now, we are concerned with a  {\it rational newform} $f$: a normalized cusp form of weight $2$ which has the rational Fourier expansion, is a simultaneous eigenform for all the Hecke operators and is a newform in the sense of  [AL].   
Let $\delta_{N}$ denote the character defined by $\delta_{N}(p)=1$ if $p\nmid N$ and $=0$ if $p\mid N$. 
\begin{prop}
Let $f(z)=\sum_{n=1}^{\infty}a_{n}(f)q^{n}$ be a rational newform. Then, the Fourier expansion of $f(z)$ satisfies the following conditions. 
\begin{enumerate}
\item $a_{p^{r+1}}(f)=a_{p}(f)a_{p^{r}}(f)-\delta_{N}(p)pa_{p^{r-1}}(f)\qquad (r\geq 1)$
\item $a_{mn}(f)=a_{m}(f)a_{n}(f)$ $\mspace{91mu}$ $\qquad ((m,n)=1)$.
\end{enumerate}
\end{prop}
Given a rational newform $f$, we consider an associated period lattice 
$$\Lambda_{f}=\verb+{+\int_{\alpha}^{\beta}f(z)dz\mid \alpha,\ \beta\in\mathbb{H}^{*}, \ \alpha\equiv \beta \ (\text{mod} \ \Gamma_{0}(N))\verb+}+$$
which is a discrete subgroup of $\mathbb{C}$ of rank $2$. Then, it is known that the quotient $E_{f}=\mathbb{C}/\Lambda_{f}$ is an elliptic curve over $\mathbb{Q}$ of conductor $N$ and that we have $L(E_{f},s)=L(f,s)$ where the LHS denotes the Hasse-Weil L-function of $E_{f}$ and the RHS denotes the Dirichlet L-series of $f$. Conversely, for any elliptic curve $E$ over $\mathbb{Q}$, there exists a rational newform $f$ such that we have $L(E,s)=L(f,s)$ ([Wi], [TW], [BCDT]). From this equality, we have the following result. 
\begin{prop}
For any prime $p\nmid N$, we have $a_{p}(f)=1+p-\# E_{f}(\mathbb{F}_{p})$ and there exists $0\leq \theta_{p}\leq \pi$ such that $a_{p}(f)=2p^{\frac{1}{2}}  \cos (\theta_{p})$.
\end{prop}
\section{Deformation of the Fourier expansion}
In this section, we shall deform the Fourier expansion of a rational newform with respect to the arguments $\verb+{+\theta_{p}\verb+}+_{p}$ (Proposition $2.2$).
\begin{defn} Let $F(z)=\sum_{n=1}^{\infty}a_{n}(F)q^{n}$ be a formal power series in $\mathbb{C}[[q]]$ which satisfies the following conditions.
\begin{enumerate}
\item If  there exists a rational newform $f(z)$ such that we have $a_{p}(f)=a_{p}(F)$ for almost all primes $p$, put $F(z)=f(z)$. The coefficients of $F(z)$ are determined by Proposition $2.1$ and $2.2$.
\item If there does not exist such a rational newform, assume that $F(z)$ is normalized (i.e. $a_{1}(F)=1$) and that, for each prime $p$, there exists  $0\leq  \theta_{p}^{F} \leq \pi$ such that we have $$a_{p}(F)=2p^{\frac{1}{2}} \cos (\theta_{p}^{F}).$$ Furthermore,  the following compatible conditions are satisfied.
\begin{enumerate}
\item $a_{p^{r+1}}(F)=a_{p}(F)a_{p^{r}}(F)-pa_{p^{r-1}}(F)\qquad (r\geq 1)$
\item $a_{mn}(F)=a_{m}(F)a_{n}(F)$ $\mspace{52mu}$ $\qquad ((m,n)=1)$.
\end{enumerate}
\end{enumerate}
\end{defn}
Fix a power series $F(z)\in \mathbb{C}[[q]]$ as above. Let $\verb+{+\gamma_{i}\verb+}+_{i=1, 2}$  denote any smooth path from $\alpha_{i}$ to
$\beta_{i}$ in $\mathbb{H}^{*}$. Consider an associated period lattice
$$\Lambda_{F}(\gamma_{1},\gamma_{2})=\verb+{+\int_{\alpha_{i}}^{\beta_{i}}F(z)dz\mid \alpha_{i}\overset{\mathrm{\gamma_{i}}}{\thicksim }\beta_{i}\verb+}+_{i=1,2}.$$ 
Note that, contrary to $\Lambda_{f}$, this  $\Lambda_{F}(\gamma_{1},\gamma_{2})$ does not form a discrete subgroup of $\mathbb{C}$ depending on the choice of $\verb+{+\gamma_{i}\verb+}+_{i=1, 2}$. Thus, the quotient $E_{F}(\gamma_{1},\gamma_{2})=\mathbb{C}/\Lambda_{F}(\gamma_{1},\gamma_{2})$ is not an elliptic curve in general.
\begin{defn}With notation as above, let $\Theta$ denote the topological family $\verb+{+E_{F}(\gamma_{1},\gamma_{2})\verb+}+$ where $F$ (resp. $\verb+{+\gamma_{i}\verb+}+_{i=1,2}$) runs through any power series as in Definition $3.1$ (resp. any smooth path in $\mathbb{H}^{*}$).
\end{defn}
\begin{rem}
We can say that this topological family $\Theta$ is the smallest in the sense that it contains all elliptic curves over $\mathbb{Q}$ up to isogeny  and the associated rational newforms are all parametrized by the arguments $\verb+{+\theta_{p}\verb+}+_{p}$.
\end{rem}
\section{Applications}
\subsection{The case of rank $0$}
For any elliptic curve $E$ over $\mathbb{Q}$, the Birch and Swinnerton-Dyer conjecture predicts that the rank of Mordell-Weil group $E(\mathbb{Q})$ is equal to the order of the zero of $L(E,s)$ at $s=1$. In the case that we have $L(E,1)\not=0$, it is known that the Mordell-Weil  group of $E$ is of rank $0$ ([CW]). Now, assume that $E$ is such an elliptic curve and that $f$ is an associated rational newform satisfying $L(E,s)=L(f,s)$. Since the Dirichlet L-series $L(f,s)$ can be written via Mellin transform
$$L(f,s)=(2\pi)^{s}\Gamma(s)^{-1}\int_{0}^{i\infty}(-iz)^{s}f(z)\frac{dz}{z}$$
where $\Gamma(s)$ denotes the gamma function of $s$, the period integral $\int_{0}^{i\infty}f(z)dz$ does not vanish. Let  $I$ denote any smooth path from $0$ to $i\infty$ in $\mathbb{H}^{*}$.
\begin{ex}
Let $\verb+{+E_{i}\verb+}+_{i=1,2}$ be two elliptic curves over $\mathbb{Q}$. Assume that  there exist a set of formal power series $\verb+{+F(z)\verb+}+_{F}$ as in Definition $3.1$ and a set of smooth paths $\verb+{+J\verb+}+_{J}$ in $\mathbb{H}^{*}$ such that $\verb+{+E_{F}(I,J)\verb+}+_{F,J}$ forms a topological family of (non-degenerate) elliptic curves connecting $E_{1}$ and $E_{2}$. Then, Mordell-Weil groups of  $\verb+{+E_{i}\verb+}+_{i=1,2}$ are of rank $0$.
\end{ex}  
\subsection{The case of rank $1$}
First, we shall recall the results of [GZ]. 
Let $K$ be an imaginary quadratic field whose discriminant $D$ is relatively prime to the level $N$ of the rational newform $f$ and let $H$ denote the Hilbert class field of $K$. 
Fix an element $\sigma$ in $\Gal(H/K)$. Note that this Galois group is isomorphic to the class group $\Cl_{K}$ of $K$. Let $\mathscr{A}_{K}$ be the class corresponding to $\sigma$ and let $\theta_{\mathscr{A}_{K}}(z)$  denote the theta series
$$\theta_{\mathscr{A}_{K}}(z)=\sum_{n\geq 0}r_{\mathscr{A}_{K}}(n)q^{n} \quad (q=e^{2\pi iz})$$
where $r_{\mathscr{A}_{K}}(0)=\frac{1}{\sharp (\mathscr{O}_{K}^{*})}$ ($\mathscr{O}_{K}$: the ring of integers in $K$) and  $r_{\mathscr{A}_{K}}(n)$ ($n\geq 1$) is the number of integral ideals $\alpha$ in the class of $\mathscr{A}_{K}$ with norm $n$. 
Define the $L$-function associated to the rational newform $f=\sum_{n}a_{n}q^{n}\in S_{2}(N)$ and the ideal class $\mathscr{A}_{K}$ by
$$L_{\mathscr{A}_{K}}(f,s)=\biggl(\sum_{n\geq 1, (n,DN)=1}\epsilon_{K}(n)n^{1-2s}\biggr)\cdot\biggl(\sum_{n\geq 1}a_{n}r_{\mathscr{A}_{K}}(n)n^{-s}\biggr)$$
where $\epsilon_{K}:(\mathbb{Z}/D\mathbb{Z})^{*}\rightarrow \verb+{+\pm 1\verb+}+$ denotes the character associated to $K/\mathbb{Q}$. Furthermore, for a complex character $\chi$ of the ideal class group of $K$, denote the total $L$-function by 
$$L(f, \chi, s)=\sum_{\mathscr{A}_{K}}\chi(\mathscr{A}_{K})L_{\mathscr{A}_{K}}(f,s).$$
Then, it is known that both of $L_{\mathscr{A}_{K}}(f,s)$ and $L(f, \chi, s)$ have analytic continuations to the entire plane and satisfy functional equations ($s \leftrightarrow 2-s$).  Furthermore, if we put $L_{\epsilon_{K}}(f,s)=\sum_{n}\epsilon_{K}(n)a_{n}n^{-s}$ for $f=\sum_{n}a_{n}q^{n}$, we have $L(f,s)L_{\epsilon_{K}}(f,s)=L(f,\text{\bf{1}},s)$. Note that  $L_{\epsilon_{K}}(f,s)$ is the Hasse-Weil $L$-function of  $E'$  over $\mathbb{Q}$ where $E'$ denotes the twist of $E$ over $K$ ([GZ, p.309, 312]). The following thing is one of the main results of Gross-Zagier.
\begin{prop}{\text{([GZ, p.230])}}
There exists a cusp form $g_{\mathscr{A}_{K}}$ of weight $2$ on $\Gamma_{0}(N)$ such that we have
$$L_{\mathscr{A}_{K}}'(f,1)=32\pi^{2}\sharp(\mathscr{O}_{K}^{*})^{-2}\verb+|+D\verb+|+^{-\frac{1}{2}}\cdot(g_{\mathscr{A}_{K}},f)_{N}$$
where $(\ ,\ )_{N}$ denotes the Petersson inner product on cusp forms of weight $2$ for $\Gamma_{0}(N)$. Thus, this formula leads to
$$L'(f,\chi,1)=\sum_{\mathscr{A}_{K}}\chi(\mathscr{A}_{K})L_{\mathscr{A}_{K}}'(f,1)=32\pi^{2}\sharp(\mathscr{O}_{K}^{*})^{-2}\verb+|+D\verb+|+^{-\frac{1}{2}}\cdot(\sum_{\mathscr{A}_{K}}\chi(\mathscr{A}_{K})g_{\mathscr{A}_{K}},f)_{N}$$
\end{prop} 
Now, let $E$ be an elliptic curve over $\mathbb{Q}$ such that  $L(E,s)=L(f,s)$ for some rational newform $f\in S_{2}(N)$. Assume that we have 
$\ord_{s=1}L(E,s)=1$. In this case, it is known that the Mordell-Weil group of $E$ is of rank $1$ ([Ko]). Furthermore, since the sign of the functional equation of $L(E,s)=L(f,s)$ is $-1$, we can choose an imaginary quadratic extension  $K/\mathbb{Q}$ such that $L_{\epsilon_{K}}(f,1)\not=0$ ([Wa]).  In particular, it follows that we obtain $L'(f,\text{\bf{1}},1)\not=0$ and thus $(\sum_{\mathscr{A}_{K}}\text{\bf{1}}(\mathscr{A}_{K})g_{\mathscr{A}_{K}},f)_{N}\not=0$. Let $\verb+{+g_{i}\verb+}+_{i=1}^{d}$ (resp. $\verb+{+h_{j}\verb+}+_{j=1}^{e}$) denote a basis of the space of newforms (resp. oldforms) in $S_{2}(N)$ over $\mathbb{C}$. If we write
$\sum_{\mathscr{A}_{K}}\text{\bf{1}}(\mathscr{A}_{K})g_{\mathscr{A}_{K}}=\sum_{i=1}^{d}a_{i}g_{i}+\sum_{j=1}^{e}b_{j}h_{j}$ $(a_{i}, b_{j}\in\mathbb{C})$, put  $G_{K}=\sum_{i=1}^{d}a_{i}g_{i}\in S_{2}(N)$.
\begin{defn}
Let $F(z)\in \mathbb{C}[[q]]$ ($q=e^{2\pi iz}$) be a formal power series  as in Definition $3.1$. Fix a fundamental domain $R$ in $\mathbb{H}$ for $\Gamma_{0}(N)$. We say that $F(z)$ is of level $N$ with respect to $R$ if we have
$$(G_{K},F(z))_{N,R}:=\displaystyle\int_{R} G_{K}\cdot\overline{F(z)}dxdy\not=0\quad (z=x+iy)$$
for some imaginary quadratic extension $K/\mathbb{Q}$ whose discriminant is relatively prime to $N$.
\end{defn}
\begin{ex}  
Let us consider the following two cases. 
\begin{enumerate} 
\item
Let  $\verb+{+F(z)\verb+}+_{F}$  be a set of formal power series  of level $N$ with respect to $R$ such that we have $L(F,1):=-2\pi i\Gamma(1)^{-1}\int_{0}^{i\infty}F(z)dz=0$ and let  $\verb+{+I,J\verb+}+_{I,J}$ denote a set of smooth paths in $\mathbb{H}^{*}$.
Assume that two  elliptic curves $\verb+{+E_{i}\verb+}+_{i=1,2}$  over $\mathbb{Q}$ of conductor $N$ are connected by the topological family $\verb+{+E_{F}(I,J)\verb+}+_{F,I,J}$.  Then, Mordell-Weil groups of  $\verb+{+E_{i}\verb+}+_{i=1,2}$ are of rank $1$. 
\item 
On the other hand, let $\mathbb{E}_{1}$ (resp. $\mathbb{E}_{2}$) be an elliptic curve over $\mathbb{Q}$ of conductor $N$ (resp. $N'$).  Here, $N'$ denotes a positive integer such that $N'\verb+|+N$ and $N'<N$. Assume that the Mordell-Weil group of $\mathbb{E}_{1}$ is of rank $1$. Then, though it may happen that the Mordell-Weil group of $\mathbb{E}_{2}$ is also of rank $1$, there is not a set of formal power series  of level $N$ connecting both elliptic curves.
\end{enumerate}
 In fancy language, we can say that the existence of  (non-torsion)  rational points on elliptic curves is partially governed  by the {\it singular locus} of special fibers in $\Spec(\mathbb{Z})$. 
\end{ex} 
\begin{rem}
Let $\verb+{+E_{i}\verb+}+_{i=1,2}$ be two elliptic curves over $\mathbb{Q}$ of conductor $N$ whose Mordell-Weil groups are of rank $1$. Take 
rational newforms $\verb+{+f_{i}\verb+}+_{i=1,2}\in S_{2}(N)$ such that we have $L(f_{i},s)=L(E_{i},s)$. Assume that the strong Birch and 
Swinnerton-Dyer conjecture holds ([C]). From the equality $L'(f_{i},1)L_{\epsilon_{K_{i}}}(f_{i},1)=L'(f_{i},\text{\bf{1}},1)$, we obtain 
$L'(f_{i},\text{\bf{1}},1)>0$ and thus $(G_{K_{i}},f_{i})_{N,R}>0$. Here, we choose imaginary quadratic fields $K_{i}/\mathbb{Q}$ such that we have  $L_{\epsilon_{K_{i}}}(f_{i},1)\not=0$. Define a set of formal power series by 
$$F_{t}(z)=tf_{1}(z)+(1-t)f_{2}(z)\quad (0\leq t\leq 1).$$
If we can take $K_{1}=K_{2}$ (e.g. two elliptic curves of conductor $91$ and $\mathbb{Q}(\sqrt{-3})$ [C, p.118 and 223-224]), we obtain  $(G_{K_{i}},F_{t}(z))_{N,R}>0$ for all $0\leq t\leq 1$. Thus,  though  this set of formal power series  $\verb+{+F_{t}(z)\verb+}+_{0\leq t\leq 1}$ (regrettably) does not satisfy the compatible conditions in Definition $3.1$, two elliptic curves $\verb+{+E_{i}\verb+}+_{i=1,2}$ are connected by this set of formal power series {\it of level $N$} anyway.  
\end{rem}

\end{document}